\nonstopmode
\documentclass[a4paper, 10pt]{amsart}
\usepackage{latexsym}
\usepackage{fancyhdr}
\usepackage{amsmath, amssymb}
\usepackage[ansinew]{inputenc}
\usepackage[all]{xy}
\usepackage{verbatim}
\usepackage{psfrag}
\theoremstyle{plain}
\newtheorem*{lemma*}{Lemma}
\newtheorem{lemma}[subsection]{Lemma}
\newtheorem*{theorem*}{Theorem}

\newtheorem*{proposition*}{Proposition}

\newtheorem*{corollary*}{Corollary}

\theoremstyle{definition}
\newtheorem*{definition*}{Definition}

\newtheorem*{example*}{Example}

\newtheorem*{remark*}{Remark}
\newtheorem*{remarks*}{Remarks}

\numberwithin{equation}{subsection}

\def\al{\alpha}
\def\be{\beta}
\def\ga{\gamma}
\def\de{\delta}
\def\ep{\epsilon}

\def\ze{\zeta}

\def\th{\theta}

\def\rh{\rho}
\def\vr{\varrho}
\def\si{\sigma}

\def\ta{\tau}

\def\vh{\varphi}
\def\ch{\chi}
\def\ps{\psi}
\def\om{\omega}

\def\De{\Delta}
\def\Th{\Theta}

\def\Si{\Sigma}

\def\Ph{\Phi}

\def\Om{\Omega}

\def\C{\mathbb{C}}

\def\N{\mathbb{N}}

\def\R{\mathbb{R}}

\def\cF{\mathcal{F}}

\def\cH{\mathcal{H}}

\def\cL{\mathcal{L}}

\def\p{\partial}

\def\<{\langle}
\def\>{\rangle}
\renewcommand{\o}{\circ}

\newenvironment{demo}[1]{\par\smallskip\noindent{\bf #1.}}{\par\smallskip}

\let\on=\operatorname

\providecommand{\bysame}{\leavevmode\hbox to3em{\hrulefill}\thinspace}
\providecommand{\MR}{\relax\ifhmode\unskip\space\fi MR }

\providecommand{\href}[2]{#2}

\title[Invariant functions in Denjoy--Carleman classes]
{Invariant functions in Denjoy--Carleman classes}

\author[A. Rainer]
{Armin Rainer}

\address{Armin Rainer: Dipartimento di Matematica, Universit\'a di Pisa, 
Largo Bruno Pontecorvo~5, 56127 Pisa, Italy}

\email{armin.rainer@univie.ac.at}

\begin{document}

\begin{abstract}
Let $V$ be a real finite dimensional representation of a compact Lie group $G$. 
It is well-known that the algebra $\R[V]^G$ of $G$-invariant polynomials on $V$ is finitely generated,
say by $\si_1,\ldots,\si_p$. Schwarz \cite{Schwarz75} proved that each $G$-invariant $C^\infty$-function 
$f$ on $V$ has the form $f=F(\si_1,\ldots,\si_p)$ for a $C^\infty$-function $F$ on $\R^p$. 
We investigate this representation within the framework of Denjoy--Carleman classes.
One can in general not expect that $f$ and $F$ lie in the same Denjoy--Carleman class $C^M$ (with $M=(M_k)$).  
For finite groups $G$ and (more generally) for polar representations $V$ we show that for each $G$-invariant 
$f$ of class $C^M$ there is an $F$ of class $C^N$ such that $f=F(\si_1,\ldots,\si_p)$, 
if $N$ is strongly regular and satisfies 
\[
 \sup_{k \in \N_{>0}} \Big(\frac{M_{km}}{N_k}\Big)^{\frac{1}{k}} < \infty,
\]
where $m$ is an (explicitly known) 
integer depending only on the representation.  
In particular, each $G$-invariant $(1+\de)$-Gevrey function $f$ (with $\de>0$) has the form $f=F(\si_1,\ldots,\si_p)$
for a $(1+\de m)$-Gevrey function $F$.
Applications to equivariant functions and basic differential forms are given.   
\end{abstract}

\thanks{The author was supported by
`Fonds zur F\"orderung der wissenschaftlichen Forschung, Projekt P19392 \& Projekt J2771'}
\keywords{Denjoy--Carleman class, invariant functions}
\subjclass[2000]{26E10, 58C25, 57S15}
\date{May 8, 2008}

\maketitle

\section{Introduction}

Let $V$ be a real finite dimensional representation of a compact Lie group $G$.
By a classical theorem due to Hilbert the algebra $\R[V]^G$ of $G$-invariant polynomials on $V$ is finitely generated. Choose a system of homogeneous generators $\si_1,\ldots,\si_p$ of $\R[V]^G$ and define $\si:=(\si_1,\ldots,\si_p) : V \to \R^p$.
Schwarz \cite{Schwarz75} proved a smooth analog of Hilbert's theorem for orthogonal representations $V$ 
of compact Lie groups $G$:
the induced mapping $\si^* : C^\infty(\R^p) \to C^\infty(V)^G$ is surjective. 
Mather \cite{Mather77} showed that this mapping is even split surjective. 

The finitely differentiable case was studied, too: $\si^* : C^n(\R^p) \to C^n(V)^G$ is in general not 
surjective, but $\si^* C^n(\R^p)$ contains $C^{n q}(V)^G$ for a suitable integer $q$. 
See \cite{Barbancon72}, \cite{BR83}, \cite{Barbancon86}, \cite{Rumberger98}. 

In this paper we treat Schwarz's theorem in the framework of Denjoy--Carleman classes. 
These classes of smooth functions play an important role in harmonic analysis and various branches of 
differential equations (especially Gevrey classes).
Let $M=(M_k)_{k \in \N}$ be a non-decreasing sequence of real numbers with $M_0=1$.
A smooth function $f$ in an open subset $U \subseteq \R^n$ belongs to the Denjoy--Carleman class $C^M(U)$ 
if for any compact subset $K \subseteq U$ there exist positive constants $C$ and $\vr$ such that
\[
|\p^\al f(x)| \le C \vr^{|\al|} |\al|! M_{|\al|}
\]
for all $\al \in \N^n$ and $x \in K$. See section \ref{secDC} for more on Denjoy--Carleman classes.
As examples (\cite{Bronshtein87}, see also \ref{sym}) show, one cannot expect in general that a smooth $G$-invariant 
function $f$ on $V$ of class $C^M$ has the form $f = F \o \si$ for a function $F$ of the same class $C^M$.
 
For finite groups $G$ and (more generally) for polar representations $V$ we prove that the representation 
$f = F \o \si$ holds in the context of Denjoy--Carleman classes, where $F$ has lower regularity than $f$. 
More precisely: Let $G$ be a subgroup of finite order $m$ of $\on{GL}(V)$. Let $M$ and $N$ be sequences 
satisfying some mild conditions which guarantee stability under composition and derivation for $C^M$ and $C^N$ (see \ref{2.1}). 
Assume that $N$ is strongly regular (see \ref{streg}) and that 
\[
 \sup_{k \in \N_{>0}} \Big(\frac{M_{km}}{N_k}\Big)^{\frac{1}{k}} < \infty.
\]
Then for any $G$-invariant function $f \in C^M(V)$ there exists a function $F \in C^N(\R^p)$ 
such that $f = F \o \si$. 
In particular: Any $G$-invariant Gevrey function $f \in G^{1+\de}(V)$ (with $\de>0$) has the form $f = F \o \si$ 
for a Gevrey function $F \in G^{1+\de m}(\R^p)$. See theorem \ref{ifDC}.
The result does not depend on the choice of generators $\si_i$, since any two choices differ only by 
a polynomial diffeomorphism and the involved Denjoy--Carleman classes are stable under composition.

Note that Thilliez \cite{Thilliez97} treats a very similar problem: For a compact subset $E \subseteq \R^n$,  an analytic mapping $\Ph : U \to \R^n$ on an open neighborhood $U$ of $E$, and a function $f \in C^M(U)$ of the form $f = g \o \Ph$ with $g \in C^\infty(W)$ for an open neighborhood $W$ of $\Ph(E)$, the existence of a sequence $N$ such that $g \in C^N(W)$ is investigated. This is done by studying the complex setting: Now $E$ is compact in $\C^n$, $\Ph$ is a $\C^n$-valued holomorphic mapping defined near $E$, and $g$ is $C^\infty$ on $\C^n$ and $\bar \p$-flat on $\Ph(E)$. 
However, our results are not covered by Thilliez', since the minimal number of generators of $\R[V]^G$ does in general not coincide with the dimension of the representation space $V$. 

We prove the main theorem in section \ref{secif}. We shall deduce it from an analog theorem (see \ref{sym}) due to 
Bronshtein \cite{Bronshtein86, Bronshtein87} which treats the standard representation of the symmetric group $\on{S}_n$ in $\R^n$.
This method is inspired by Barbançon and Raïs \cite{BR83} 
deploying Weyl's account \cite{Weyl39} of Noether's \cite{Noether16} proof of Hilbert's theorem.

The rest of the paper is devoted to several applications of this main theorem. In section \ref{secequi} 
we treat the presentation in Denjoy--Carleman classes of equivariant functions between representations of a finite group. 

In section \ref{secpolar} the main theorem \ref{ifDC} is generalized to polar representations, i.e., orthogonal 
finite dimensional representations $V$ of a compact Lie group $G$ allowing a linear subspace $\Si \subseteq V$ 
which meets each orbit orthogonally (see theorem \ref{polDC}). 
The trace of the $G$-action in $\Si$ is the action of the generalized Weyl group $W$ which is a finite group.
In analogy with a result due to Palais and Terng \cite{PT87}, which states that restriction induces an isomorphism $I_1 : C^\infty(V)^G \cong C^\infty(\Si)^W$, we show that each $W$-invariant function on $\Si$ of class 
$C^M$ has a $G$-invariant extension to $V$ of class $C^N$, where $M$ and $N$ are sequences with the aforementioned
properties. More generally, Michor \cite{Michor96B, Michor97B} proved that restriction induces an isomorphism 
$I_2 : \Om_{\on{hor}}^p(V)^G \cong \Om^p(\Si)^W$, where $\Om_{\on{hor}}^p(V)^G$ is the space of basic $p$-forms on $V$, 
i.e., $G$-invariant forms that kill each vector tangent to some orbit. 
Our main theorem \ref{ifDC} allows to conclude that each $W$-invariant $p$-form on $\Si$ of class 
$C^M$ has a basic extension to $V$ of class $C^N$ (with $M$ and $N$ as above).

In \cite{PT87} and  \cite{Michor96B, Michor97B} the isomorphisms $I_1$ and $I_2$ are established in the more general setting of smooth proper Riemannian $G$-manifolds $X$ with sections, where there exist closed submanifolds $\Si \subseteq X$ meeting each orbit orthogonally.
In section \ref{secGmf} we explain that our analog results in the framework of Denjoy--Carleman classes generalize 
to real analytic proper Riemannian $G$-manifolds $X$ with sections.

\section{Denjoy--Carleman classes} \label{secDC}

\subsection{\label{2.1}Denjoy--Carleman classes of differentiable functions} 

We mainly follow \cite{Thilliez08} (see also the references therein).
We use $\N = \N_{>0} \cup \{0\}$.
For each multi-index $\al=(\al_1,\ldots,\al_n) \in \N^n$, we write
$\al!=\al_1! \cdots \al_n!$, $|\al|= \al_1 +\cdots+ \al_n$, and 
$\p^\al=\p^{|\al|}/\p x_1^{\al_1} \cdots \p x_n^{\al_n}$. 

Let $M=(M_k)_{k \in \N}$ be an increasing sequence ($M_{k+1}\ge M_k$) 
of real numbers with $M_0=1$.
Let $U \subseteq \R^n$ be open.
We denote by $C^M(U)$ the set of all $f \in C^\infty(U)$ such that, for all
compact $K \subseteq U$,
there exist positive constants $C$ and $\vr$ such that
\begin{equation}\label{CM}
|\p^\al f(x)| \le C \, \vr^{|\al|} \, |\al|! \, M_{|\al|}
\end{equation}
for all $\al \in \N^n$ and $x \in K$.
The set $C^M(U)$ is the \emph{Denjoy--Carleman class} of functions on $U$.
If $M_k=1$, for all $k$, then $C^M(U)$ coincides with the ring $C^\om(U)$
of real analytic functions
on $U$. In general, $C^\om(U) \subseteq C^M(U) \subseteq C^\infty(U)$.

We assume that $M=(M_k)$ is \emph{logarithmically convex}, i.e.,
\begin{equation} \label{logconvex}
M_k^2 \le M_{k-1} \, M_{k+1} \quad \text{ for all } k,
\end{equation}
or, equivalently, $M_{k+1}/M_k$ is increasing.
Considering $M_0=1$, we obtain that also $(M_k)^{1/k}$ is increasing and
\begin{equation}
M_l \, M_k\le M_{l+k} \quad \text{  for all }l,k\in \N.
\label{logconvex1}
\end{equation}

Hypothesis \eqref{logconvex} implies that $C^M(U)$ is a ring, for all open
subsets $U \subseteq \R^n$, 
which can easily be derived from \eqref{logconvex1} by means of Leibniz's rule.
Note that definition \eqref{CM} makes sense also for functions 
$U\to \mathbb R^p$.
For $C^M$-mappings, \eqref{logconvex} guarantees stability under composition
(\cite{Roumieu62/63}, see also \cite[4.7]{BM04}).

A further consequence of \eqref{logconvex} is the inverse function theorem
for $C^M$ (\cite{Komatsu79}; for a proof see also \cite[4.10]{BM04}): Let
$f : U \to V$ be a $C^M$-mapping between open subsets $U,V \subseteq \R^n$.
Let $x_0 \in U$.
Suppose that the Jacobian matrix $(\p f/\p x)(x_0)$ is invertible. Then
there are neighborhoods $U'$ of $x_0$, $V'$ of
$y_0 := f(x_0)$ such that $f:U'\to V'$ is a $C^M$-diffeomorphism.

Moreover, \eqref{logconvex} implies that $C^M$ is closed under solving ODEs
(due to \cite{Komatsu80}).

Suppose that $M=(M_k)$ and $N=(N_k)$ satisfy $M_k \le C^k \, N_k$, for all $k$
and a constant $C$, or equivalently,
\begin{equation} \label{incl}
\sup_{k \in \N_{>0}} \Big(\frac{M_k}{N_k}\Big)^{\frac{1}{k}} < \infty.
\end{equation}
Then, evidently $C^M(U) \subseteq C^N(U)$. The converse
is true as well (if \eqref{logconvex} is assumed): 
One can prove that there exists $f \in C^M(\R)$ such that 
$|f^{(k)}(0)| \ge k! \, M_k$ for all $k$ (see \cite[Theorem 1]{Thilliez08}).
So the inclusion $C^M(U) \subseteq C^N(U)$ implies \eqref{incl}. 

Setting $N_k=1$ in \eqref{incl} yields that $C^\om(U) = C^M(U)$ if and only if 
\[
\sup_{k \in \N_{>0}} (M_k)^{\frac{1}{k}} < \infty.
\] 
Since $(M_k)^{1/k}$ is increasing (by logarithmic convexity), 
the strict inclusion $C^\om(U) \subsetneq C^M(U)$ is equivalent to  
\[
\lim_{k \to \infty} (M_k)^{\frac{1}{k}} = \infty.
\] 

We shall also assume that $C^M$ is stable under derivation, which is
equivalent to the following condition
\begin{equation} \label{der}
\sup_{k \in \N_{>0}} \Big(\frac{M_{k+1}}{M_k}\Big)^{\frac{1}{k}} < \infty.
\end{equation}
Note that the first order partial derivatives of elements in $C^M(U)$
belong to $C^{M^{+1}}(U)$,
where $M^{+1}$ denotes the shifted sequence 
$M^{+1} = (M_{k+1})_{k \in \N}$.
So the equivalence follows from \eqref{incl}, by replacing $M$ with
$M^{+1}$ and $N$ with $M$.

\subsection*{Definition} By a {\it DC-weight sequence} we mean a sequence 
$M=(M_k)_{k \in \N}$ of positive numbers with $M_0=1$ 
which is monotone increasing
($M_{k+1}\ge M_k$), logarithmically convex
\eqref{logconvex}, and satisfies \eqref{der}.
Then $C^M(U,\mathbb R)$ is a differential ring, and the class of
$C^M$-functions is stable under compositions, as above.

\subsection{\label{2.2}Quasianalytic function classes}
Let $\cF_n$ denote the ring of formal power series in $n$ variables (with
real or complex coefficients).
We denote by $\cF_n^M$ the set of elements 
$F=\sum_{\al \in \N^n} F_\al \, x^\al$ of $\cF_n$ for which there 
exist positive constants $C$ and $\vr$ such that
\[
|F_\al| \le C \, \vr^{|\al|} \, M_{|\al|}
\]
for all $\al \in \N^n$.
A class $C^M$ is called \emph{quasianalytic} if, for open connected $U
\subseteq \R^n$ and all $a \in U$, the Taylor series homomorphism
\[
T_a : C^M(U) \to \cF_n^M,
~f \mapsto T_a f(x) = \sum_{\al \in \N^n} \frac{1}{\al!} \, \p^\al f(a) \, x^\al
\]
is injective.
By the Denjoy--Carleman theorem (\cite{Denjoy21}, \cite{Carleman26}),
$C^M$ is quasianalytic if and only if
\begin{equation} \label{eqqa}
\sum_{k=0}^\infty \frac{M_k}{(k+1) \, M_{k+1}}=\infty,
\quad\text{  or, equivalently, }\quad
\sum_{k=1}^\infty \Big(\frac1{k!\,M_k}\Big)^{\frac1k} = \infty.
\end{equation}
For contemporary proofs see for instance \cite[1.3.8]{Hoermander83I} 
or \cite[19.11]{Rudin87}. 

Suppose that $C^\om(U) \subsetneq C^M(U)$ and $C^M(U)$ is quasianalytic. 
Then $T_a : C^M(U) \to \cF_n^M$ is not surjective. 
This is due to Carleman \cite{Carleman26}; 
an elementary proof can be found in \cite[Theorem 3]{Thilliez08}.

\subsection{\label{2.3}Non-quasianalytic function classes}
If $M$ is a DC-weight sequence which is not quasianalytic, then there are 
$C^M$ partitions of unity. Namely, there exists a $C^M$ function $f$ on $\mathbb R$ 
which does not vanish in any neighborhood of 0 but which has vanishing
Taylor series at 0. Let $g(t)=0$ for $t\le 0$ and $g(t)=f(t)$ for $t>0$.
 From $g$ we can construct $C^M$ bump functions as usual. 

\subsection{\label{2.4}Strong non-quasianalytic function classes}
Let $M$ be a DC-weight sequence with $C^\om(U, \mathbb R) \subsetneq C^M(U,\mathbb R)$. 
Then the mapping $T_a : C^M(U,\mathbb R) \to \cF_n^M$
is surjective, for all $a \in U$,
if and only if there is a constant $C$ such that
\begin{equation} \label{strnonqa}
\sum_{k=j}^\infty \frac{M_k}{(k+1) \, M_{k+1}} \le C \frac{M_j}{M_{j+1}} \quad
\text{for any integer}~ j \ge 0.
\end{equation}
See \cite{Petzsche88} and references therein.
\eqref{strnonqa} is called \emph{strong non-quasianalyticity} condition.

\subsection{\label{2.5}Moderate growth}
A DC-weight sequence $M$ has
\emph{moderate growth} if
\begin{equation} \label{mgrowth}
\sup_{j,k \in \N_{>0}} \Big(\frac{M_{j+k}}{M_j \, M_k}\Big)^{\frac{1}{j+k}} <
\infty.
\end{equation}

\subsection{\label{streg}Strong regularity}
Moderate growth \eqref{mgrowth} together with 
strong non-quasi\-analyticity \eqref{strnonqa} is called
\emph{strong regularity}: Then
a version of Whitney's
extension theorem holds for
the corresponding function classes. 

\subsection{\label{CMjets}Whitney's extension theorem}
Let $K \subseteq \R^n$ be compact. Denote by $J^\infty(K)$ the $C^\infty$ Whitney jets on $K$. 
We say that $F =(F_\al)_{\al \in \N^n} \in J^\infty(K)$ is a \emph{$C^M$-jet} on $K$, or belongs to $J^M(K)$, if
there exist positive constants $C$ and $\vr$ such that 
\begin{equation} \label{jet1}
 |F_\al(x)| \le C \vr^{|\al|} |\al|! \, M_{|\al|}
\end{equation}
for all $\al \in \N^n$ and $x \in K$ and 
\begin{equation} \label{jet2}
 |F_\be(x) - \p^\be T_a^p F(x)| \le C \vr^p |\be|!\, M_{p+1} |x-a|^{p+1-|\be|}
\end{equation}
for all $p \in \N$, all $\be \in \N^n$ with $|\be| \le p$ and all $x \in K$,
where
\[
 T_a^p F(x) = \sum_{|\be| \le p} \frac{1}{\be!} F_\be(a) (x-a)^\be.
\]
If $M$ is strongly regular then a version of Whitney's extension theorem holds (see \cite{Bruna80},  \cite{BBMT91}, and \cite{ChaumatChollet92}):
the mapping $J_K : C^M(\R^n) \to J^M(K), f \mapsto (\p^\al f|_K)_{\al \in \N^n}$ is surjective.

Note that, if $f \in C^\infty(\R^n)$ such that $F=J_K f$ satisfies \eqref{jet1} and if $K$ is Whitney $1$-regular, then \eqref{jet2} is 
automatically fulfilled (see \cite[3.12]{BBMT91}). Recall that $K$ is \emph{Whitney $1$-regular} if any two points $x$ and $y$ in $K$ can be connected by a path in $K$ 
of length $\le C |x-y|$, where the constant $C$ depends only on $K$.

\subsection{\label{2.6}Gevrey functions}
Let $\de >0$ and put $M_k=(k!)^\de$, for $k \in \N$.
Then $M=(M_k)$ is strongly regular. The corresponding class $C^M$ of functions is the
\emph{Gevrey class} $G^{1+\de}$.

\subsection{\label{2.7}More examples}
Let $\de >0$ and put $M_k=(\log(k+e))^{\de \, k}$, for $k \in \N$.
Then $M=(M_k)$ is quasianalytic for $0 < \de \le 1$ and non-quasianalytic 
(but not strongly) for $\de > 1$.

Let $q > 1$ and put $M_k= q^{k^2}$, for $k \in \N$. 
The corresponding $C^M$-functions are called \emph{$q$-Gevrey regular}.
Then $M=(M_k)$ is strongly non-quasianalytic but not of moderate growth, 
thus not strongly regular.

\subsection{\label{2.8}Spaces of $C^M$-functions}
Let $U \subseteq \R^n$ be open.
For any $\vr >0$ and $K \subseteq U$ compact with smooth boundary, define
\[
C^M_\vr(K) := \{f \in C^\infty(K) : \|f\|_{\vr,K} < \infty\}
\]
with
\[
\|f\|_{\vr,K} :=
\sup \Big\{ \frac{|\p^{\al} f(x)|}{\vr^{|\al|} \, |\al|! \, M_{|\al|}} : \al \in
\N^n, x \in K\Big\}.
\]
It is easy to see that $C^M_\vr(K)$ is a Banach space.
In the description of $C^M_\vr(K)$, instead of compact $K$ with smooth
boundary, we may also use open $K\subset U$ with $\overline K$ compact in
$U$, like \cite{Thilliez08}. Or we may work with Whitney jets on compact
$K$, like \cite{Komatsu73}.

The space $C^M(U)$ carries the projective limit topology over compact $K
\subseteq U$ of the inductive limit over $\vr \in \N_{>0}$:
\[
C^M(U) = \varprojlim_{K \subseteq U} \big(\varinjlim_{\vr \in \N_{>0}}
C^M_\vr(K)\big).
\]
One can prove that, for $\vr < \vr'$, the canonical injection $C^M_\vr(K) \to C^M_{\vr'}(K)$
is a compact mapping (see \cite{Komatsu73}).
Hence $\varinjlim_\vr C^M_\vr(K)$ is a Silva space, i.e., an inductive limit of
Banach spaces such that the canonical
mappings are compact.

\subsection{\label{dense}Polynomials are dense in $C^M(U)$} 
Let $M$ be a DC-weight sequence and let $U \subseteq \R^n$ be open. It is proved in \cite[3.2]{Hoermander85} (see also \cite[3.2]{HeinrichMeise07}) that the space of entire functions $\cH(\C^n)$ is dense in $C^M(U)$. Since the polynomials are dense in $\cH(\C^n)$ and the inclusion $\cH(\C^n) \to C^M(U)$ is continuous, we obtain that the polynomials are dense in $C^M(U)$. For convenience we give a proof.

\begin{lemma*}
 Let $M$ be a DC-weight sequence and let $U \subseteq \R^n$ be open. Then $\cH(\C^n)$ is dense in $C^M(U)$.
\end{lemma*}

\begin{demo}{Proof}
Let $f \in C^M(U)$ and $K \subseteq U$ compact. Let $0 < c <1$ such that $Q=K+B_c(0) \subseteq U$, where $B_c(0)=\{x\in \R^n: |x|\le c\}$. Let $\ch \in C^\infty(U)$ with $0 \le \ch \le 1$, $\ch|_Q=1$, and compact support $Q_1 = \on{supp}(\ch) \subseteq U$. 
We define for $j \in \N_{>0}$
\[
f_j := E_j * \ch f \in \cH(\C^n), \quad \text{where} \quad E_j : \C^n \to \C, ~ z \mapsto \big(\tfrac{j}{\pi}\big)^{\frac{n}{2}} e^{-j \<z \mid z\>}.
\] 
Induction shows 
\[
 \p_{i_N} \cdots \p_{i_1}(E_j * \ch f) 
\!=\! E_j *(\ch \p_{i_N} \cdots \p_{i_1}f)+\sum_{\nu=1}^N (\p_{i_N} \cdots \p_{i_{\nu+1}}E_j)*(\p_{i_\nu}\ch)(\p_{i_{\nu-1}} \cdots \p_{i_1}f),
\]
for all $N \in \N$ and $j \in \N_{>0}$, and hence
\begin{align} \label{dense1}
|\p_{i_N} \cdots \p_{i_1}(f-f_j)| \le ~&|\p_{i_N} \cdots \p_{i_1}f-E_j *(\ch \p_{i_N} \cdots \p_{i_1}f)| \nonumber\\
&{}+\sum_{\nu=1}^N |(\p_{i_N} \cdots \p_{i_{\nu+1}}E_j)*(\p_{i_\nu}\ch)(\p_{i_{\nu-1}} \cdots \p_{i_1}f)|.
\end{align}

We have for $x \in K$ and $\al \in \N^n$
\begin{align*}
 |E_j*(\ch \p^\al f)(x) &- \p^\al f(x)| = |\int E_j(y) \big(\ch(x-y) \p^\al f(x-y)-\p^\al f(x) \big) dy| \\
&\le \int_{B_c(0)} E_j(y) |\p^\al f(x-y) -\p^\al f(x)|dy \\
&\quad+{} \int_{\R^n \setminus B_c(0)} E_j(y) \big(\ch(x-y) |\p^\al f(x-y)|+|\p^\al f(x)|\big) dy.
\end{align*}
By the generalized mean value theorem we have for $x \in K$, $y \in B_c(0)$, and $\al \in \N^n$
\begin{align*}
 |\p^\al f(x-y)-\p^\al f(x)| \le \sqrt{n}\, |y| \sup_{\substack{1 \le i \le n\\0\le t\le1}} |\p^{\al+e_i} f(x-ty)|.
\end{align*}
Choose $\vr_1>0$ such that $\|f\|_{\vr_1,Q_1}<\infty$. Then for $x \in K$, $y \in B_c(0)$, and $\al \in \N^n$
\begin{align*}
 |\p^\al f(x-y)-\p^\al f(x)| &\le \sqrt{n}\, |y| \|f\|_{\vr_1,Q_1} \vr_1^{|\al|+1} (|\al|+1)!\, M_{|\al|+1} \\
&\le 2 \sqrt{n}\, |y| \vr_1  \|f\|_{\vr_1,Q_1} (2 \vr_1 C)^{|\al|} |\al|!\, M_{|\al|}  \quad(\text{by } \eqref{der}), 
\end{align*}
where $C$ is a positive constant. 
For all $j \in \N_{>0}$ we have
\[
\int_{\R^n} |y| E_j(y) dy \le \frac{C_1}{\sqrt{j}} \quad \text{and} \quad \int_{\R^n\setminus B_c(0)} E_j(y) dy \le \frac{C_1}{\sqrt{j}}
\]
for a constant $C_1$ independent of $j$. Thus there exist positive constants $C_2$ and $\vr_2$ independent of $x$, $\al$, and $j$ such that
\begin{align} \label{dense2}
 |E_j*(\ch \p^\al f)(x) &- \p^\al f(x)| \le \frac{C_2}{\sqrt{j}} \,  \vr_2^{|\al|} |\al|!\, M_{|\al|}.
\end{align}

We have for $x \in K$
\begin{align*}
|(\p_{i_N} &\cdots \p_{i_{\nu+1}}E_j)*(\p_{i_\nu}\ch)(\p_{i_{\nu-1}} \cdots \p_{i_1}f)(x)| \\
&\le 
|Q_1|
\sup_{\substack{1 \le i \le n\\y \in U}} |\p_i \ch(y)|
\sup_{y \not\in B_c(0)} |\p_{i_{N}} \cdots \p_{i_{\nu+1}} E_j(y)|
\sup_{u \in Q_1} |\p_{i_{\nu-1}} \cdots \p_{i_1} f(u)|,
\end{align*}
where $|Q_1|$ denotes the Lebesgue measure of $Q_1$. Cauchy's inequalities imply for each $\al \in \N^n$ and $r > 0$
\[
|\p^\al E_j(y)| \le \frac{\al! }{r^{|\al|}} \sup_{z \in D_r(y)} |E_j(z)|, 
\]
where $D_r(y)=\{z \in \C^n : |z_i-y_i|\le r \text{ for all } i\}$. Choosing $r=\tfrac{c}{4 \sqrt{n}}$ we get for $y \in \R^n \setminus B_c(0)$
\[
|\p^\al E_j(y)| \le \frac{\al! }{r^{|\al|}} \Big(\frac{j}{\pi}\Big)^{\frac{n}{2}} e^{-\frac{j c^2}{2}}. 
\]
Hence with $C_3=|Q_1| \|f\|_{\vr_1,Q_1} \sup_{\substack{1 \le i \le n\\y \in U}} |\p_i \ch(y)|$ we obtain for $x \in K$
\begin{align} \label{dense3}
|(\p_{i_N} \cdots \p_{i_{\nu+1}}E_j)*&(\p_{i_\nu}\ch)(\p_{i_{\nu-1}} \cdots \p_{i_1}f)(x)| \nonumber\\
&\le 
C_3
\frac{(N-\nu)!}{r^{N-\nu}} \Big(\frac{j}{\pi}\Big)^{\frac{n}{2}} e^{-\frac{j c^2}{2}}
\vr_1^{\nu-1} (\nu-1)!\, M_{\nu-1} \nonumber \\
&\le 
C_3 \Big(\frac{j}{\pi}\Big)^{\frac{n}{2}} e^{-\frac{j c^2}{2}} 
\vr_3^{N} (N-1)!\, M_{N},
\end{align}
where $\vr_3=\max\{\tfrac{1}{r},\vr_1\}$.

It follows from \eqref{dense1}, \eqref{dense2}, and \eqref{dense3} that for $\vr_4=\max\{\vr_2,\vr_3\}$
\[
\|f-f_j\|_{\vr_4,K} \le \frac{C_2}{\sqrt{j}} + C_3 \Big(\frac{j}{\pi}\Big)^{\frac{n}{2}} e^{-\frac{j c^2}{2}}.
\]
That implies the assertion.
\qed\end{demo}

\subsection{\label{clideal}Closed ideals}  
Let $U \subseteq \R^n$ be open.
Let $\vh \in C^\om(U)$. 
Consider the principal ideal $\vh C^M(U)$ generated by $\vh$.

\begin{proposition*}
Assume that $C^M$ is stable under derivation \eqref{der}. Let $\vh$ be a linear form on $\R^n$.
Then the ideal $\vh C^M(\R^n)$ is closed in $C^M(\R^n)$. More generally, assume that $\ps=\vh_1^{p_1} \cdots \vh_l^{p_l}$
is a finite product of linear forms $\vh_i$. Then $\ps C^M(\R^n)$ is closed in $C^M(\R^n)$.
\end{proposition*}

\begin{demo}{Proof}
Let $f \in \overline{\vh C^M(\R^n)}$. Then $f|_{\vh^{-1}(0)} = 0$, since evaluation at points is continuous. As $C^M$ is stable under derivation, 
the standard integral formula (after suitable linear coordinate change) implies that $f=\vh g$ for a unique 
$g \in C^M(\R^n)$. The same reasoning shows that  $\ps C^M(\R^n)$ is closed in $C^M(\R^n)$, where $\ps = \vh_1^{p_1}$. 

For the general statement it suffices to show: 
{\it Let $\ps_1$ be a polynomial and $\ps_2$ a power of a linear form.
If $\ps_1$ and $\ps_2$ are relatively prime and both generate closed ideals in $C^M(\R^n)$, 
then $\ps_1 \ps_2 C^M(\R^n)$ is closed in $C^M(\R^n)$.}
For $f \in \overline{\ps_1 \ps_2 C^M(\R^n)}$ we find functions $g_1, g_2 \in C^M(\R^n)$ with $f = \ps_1 g_1 = \ps_2 g_2$. 
Since $\ps_1$ and $\ps_2$ are relatively prime, we have $g_1|_{\ps_2^{-1}(0)}=0$. By the standard integral formula we obtain as above 
$g_1 = \ps_2 h$ with $h \in C^M(\R^n)$. Hence the assertion.
\qed\end{demo}

\begin{remark*}
Note that for any hyperbolic polynomial $\vh$ the principal ideal $\vh C^M(\R^n)$ is closed in $\vh C^M(\R^n)$ (e.g.\ \cite[4.2]{Thilliez08}). 
This follows from the fact (due to \cite{CC04}) that Weierstrass division holds in $C^M$ for hyperbolic divisors. 
A polynomial $\vh(x',x_n)=x_n^d+\sum_{j=1}^d a_j(x') x_n^{d-j}$ with $a_j \in C^M(\R^{n-1})$ and $a_j(0)=0$, for 
$1 \le j \le d$, is called \emph{hyperbolic} if,
for each $x' \in \R^{n-1}$, all roots of $\vh(x',\cdot)$ are real.

But in general the principal ideal $\vh C^M(U)$ generated by a real analytic function $\vh$ need not be closed (see \cite{Thilliez01} and \cite[part 4]{Thilliez08}). 
Compare this with the famous results on the division of distributions due to H\"ormander \cite{Hoermander58} and Lojasiewicz \cite{Lojasiewicz58,Lojasiewicz59}.
\end{remark*}

\subsection{\label{DCmf}\!\!}
Let $M$ be a DC-weight sequence, and let $X$ be a real analytic manifold. 
We can define the space $C^M(X)$ of functions of Denjoy--Carleman class $C^M$ on $X$ by means of local coordinate 
systems, since $C^M$ contains the real analytic functions and is stable under composition. 
Similarly, we may consider the space $(\Om^M)^p(X)$ of $p$-forms of class $C^M$ on $X$.

\section{Invariant functions in Denjoy--Carleman classes} \label{secif}

Throughout this paper we consider a compact Lie group $G$ acting smoothly on a manifold $X$. 
A function $f$ on $X$ is said to be 
$G$-invariant if $f(g.x)=f(x)$ for all $g \in G$ and all $x \in X$. If $\cF$ is a set of functions on $X$, then 
$\cF^G$ denotes the subset of $G$-invariant elements in $\cF$.

\subsection{Hilbert's theorem} (e.g.\ \cite{Weyl39}) \label{Hilbert}
Let $G$ be a compact Lie group and let $V$ be a real finite dimensional $G$-module. 
Then, by a theorem due to Hilbert,  the algebra $\R[V]^G$ of $G$-invariant polynomials on $V$ is finitely generated. 
The generators can be chosen homogeneous and with positive degree.

\subsection{Schwarz's theorem} \label{Schwarz}
Suppose that the representation of $G$ in $V$ is orthogonal. 
Let $\si_1,\ldots,\si_p$ be a system of generators of $\R[V]^G$ and put $\si=(\si_1,\ldots,\si_p) : V \to \R^p$.
Schwarz \cite{Schwarz75} proved that $\si^* : C^\infty(\R^p) \to C^\infty(V)^G$ is surjective, 
which is the smooth analog of \ref{Hilbert}. Mather \cite{Mather77} showed that 
$\si^* : C^\infty(\R^p) \to C^\infty(V)^G$ is even split surjective, i.e., it allows a continuous linear section.

\subsection{Symmetric functions in Denjoy--Carleman classes} \label{sym}
In the case that the symmetric group $\on{S}_n$ acts in $\R^n$ by permuting the coordinates,
the statement of Schwarz's theorem \ref{Schwarz} is due to Glaeser \cite{Glaeser63F}. 
In that case $\si_i$ is the $i$-th elementary symmetric function, i.e., 
$\si_i(x)=\sum_{1 \le j_1 < \cdots < j_i \le n} x_{j_1} \cdots x_{j_i}$, and $\si=(\si_1,\ldots,\si_n) : \R^n \to \R^n$.

The representation of symmetric functions in Denjoy--Carleman (Gevrey) classes was treated
by Bronshtein \cite{Bronshtein86, Bronshtein87}. Since we shall need it later, we present a more general version
and we sketch a proof. Let $\prod_{j=1}^p \on{S}_n$ act in $\bigoplus_{j=1}^p \R^n$ by permuting the coordinates. 
Since $\R[\bigoplus_{j=1}^p \R^n]^{\prod_{j=1}^p \on{S}_n} \cong \bigotimes_{j=1}^p \R[\R^n]^{\on{S}_n}$,
a $\prod_{j=1}^p \on{S}_n$-invariant function $f$ on $\R^{pn}$ has the form $f = F \o \th$ with 
$\th = (\si,\ldots,\si)$.

\begin{theorem*} 
Assume that $M$ and $N$ are increasing logarithmically convex sequences with $M_0=N_0=1$.
Then for any  function $f \in C^M(\R^{pn})^{\prod_{j=1}^p \on{S}_n}$ there exists a function $F \in C^N(\th(\R^{pn}))$ 
such that $f = F \o \th$ if and only if 
\begin{equation} \label{condBr}
 \sup_{k \in \N_{>0}} \Big(\frac{M_{k n}}{N_k}\Big)^{\frac{1}{k}} < \infty.
\end{equation}
\end{theorem*}

\begin{demo}{Sketch of proof}
We indicate and adapt the main steps in Bronshtein's proof. The necessity of \eqref{condBr} is shown by considering the symmetric function $f \in C^M(\R^n)$ (for $n>2$) given by
\[
f(x) = \sum_{k=0}^\infty c_k \big(1-\prod_{j=1}^n(\rh_k x_j e^{-\rh_k^2 x_j^2})\big)^{-1},
\]
where $\rh_k = \tfrac{M_{kn+1}}{M_{kn}}$ and $c_k = \tfrac{M_{kn}}{2^k \rh_k^{kn}}$. Then $f=F\o\si$ with
\[
F(\si) = \sum_{k=0}^\infty c_k \big(1-\rh_k^n \si_n e^{-\rh_k^2 (\si_1^2-2 \si_2)}\big)^{-1},
\]
and hence 
\[
 |(\p_{\si_n})^m F(0)| = \sum_{k=0}^\infty c_k \rh_k^{mn} m! \ge c_m \rh_m^{mn} m! = \frac{m! M_{mn}}{2^m}.
\]
Since $F \in C^N$ this implies \eqref{condBr}. For $n=2$ one can find a similar example.

Without loss suppose that $f \in C^M(\R^{2n})^{\on{S}_n \times \on{S}_n}$. Instead of the elementary symmetric polynomials $\si_i$ we use the Newton polynomials 
$\nu_i(x) = \sum_{j=1}^n x_j^i$ and put $\nu=(\nu_1,\ldots,\nu_n)$ (see remark \ref{ifDC}(3)). Then we may write $f(x,y)=F(\nu(x),\nu(y))=F(u,v)$ where $u=\nu(x)$, $v=\nu(y)$, and $(x,y) \in \R^n \times \R^n$. A direct computation gives
\begin{align*}
 \p_{u_k} F(u,v) &= \frac{(-1)^{k+1}}{k} \sum_{i=1}^n \frac{\si_{n-k}(x_i') \p_{x_i} f(x,y)}{\prod_{j\ne i} (x_j-x_i)} 
= \sum_{i=1}^n \frac{g_{ki}(x,y)}{\prod_{j\ne i} (x_j-x_i)},\\
 \p_{v_k} F(u,v) &= \frac{(-1)^{k+1}}{k} \sum_{i=1}^n \frac{\si_{n-k}(y_i') \p_{y_i} f(x,y)}{\prod_{j\ne i} (y_j-y_i)} 
= \sum_{i=1}^n \frac{h_{ki}(x,y)}{\prod_{j\ne i} (y_j-y_i)},
\end{align*}
where $x_i'=(x_1,\ldots,\widehat{x_i},\ldots,x_n)$, $\si_j(x_i')$ is the elementary symmetric function of degree $j$ in $n-1$ variables ($\si_0=1$), respectively for $y$, and 
\begin{align*}
 g_{ki}(x,y) &= \frac{(-1)^{k+1}}{k} \si_{n-k}(x_i') \p_{x_i} f(x,y),\\ 
h_{ki}(x,y) &= \frac{(-1)^{k+1}}{k} \si_{n-k}(y_i') \p_{y_i} f(x,y).
\end{align*}
One shows (see \cite{Bronshtein86, Bronshtein87}) that 
\begin{equation} \label{Aop}
 \p_{u_k} F = (\prod_{j=1}^{n-1} A_j^x) g_{kn} \quad \text{and} \quad \p_{v_k} F = (\prod_{j=1}^{n-1} A_j^y) h_{kn},
\end{equation}
where the operators $A_j^x$ and $A_j^y$ are defined by
\begin{align*}
 (A_j^x h)(x,y) &= \int_0^1 [(\p_{x_j}-\p_{x_{j+1}})h](t P_{j,j+1}x + (1-t)x,y) dt,\\
(A_j^y h)(x,y) &= \int_0^1 [(\p_{y_j}-\p_{y_{j+1}})h](x,t P_{j,j+1}y + (1-t)y) dt,
\end{align*}
with $P_{j,j+1}$ the linear operator $\R^n \to \R^n$ which interchanges the $j$-th and the $(j+1)$-st coordinate.

We consider
\[
 L_x^\al = \prod_{i=1}^n \p_{x_i}^{\al_i} \prod_{1 \le p < q \le n} (\p_{x_p}-\p_{x_q})^{\al_{pq}}, \quad \al=(\al_1,\al_2) \in \N^n \times \N^{\binom{n}{2}},
\]
and likewise $L_y^\al$. Let $K,L \subseteq \R^n$ be convex, compact, and $\on{S}_n$-invariant. For non-negative $m$ and $\mu$ we write
\[
 \|f\|_{\vr,K\times L}^{m,\mu} 
= \sup_{\substack{\al,\be\\(x,y) \in K \times L}} 
\frac{|L_x^\al L_y^\be f(x,y)|}{\vr^{|\al|+|\be|+m} |\al_1|! |\be_1|!\al_2!\be_2! (|\al|+|\be|+m+1)^\mu M_{|\al|+|\be|+m}}.
\]
If $f \in C^M(\R^{2n})$ then $\|f\|_{\vr,K\times L}^{m,\mu}<\infty$ for sufficiently large $\vr$. We have the following estimates
\begin{gather}
\label{pxf}
\|\p_{x_i} f\|_{\vr,K\times L}^{m+1,\mu+1} \le \|f\|_{\vr,K\times L}^{m,\mu} \quad \text{and} \quad 
\|\p_{y_i} f\|_{\vr,K\times L}^{m+1,\mu+1} \le \|f\|_{\vr,K\times L}^{m,\mu},\\
\label{xf}
\|x_i f\|_{\vr,K\times L}^{m,\mu} \le C \|f\|_{\vr,K\times L}^{m,\mu} \quad \text{and} \quad 
\|y_i f\|_{\vr,K\times L}^{m,\mu} \le C \|f\|_{\vr,K\times L}^{m,\mu},\\
\label{Af}
\|A_j^x f\|_{\vr,K\times L}^{m+1,\mu} \le C \|f\|_{\vr,K\times L}^{m,\mu} \quad \text{and} \quad 
\|A_j^y f\|_{\vr,K\times L}^{m+1,\mu} \le C \|f\|_{\vr,K\times L}^{m,\mu}.
\end{gather}
It is easy to verify \eqref{pxf} and \eqref{xf}. For the proof of \eqref{Af} we refer to \cite{Bronshtein86, Bronshtein87}. 

It follows from \eqref{Aop} and from \eqref{pxf}, \eqref{xf}, and \eqref{Af} that
\[
 \|\p_u^\al \p_v^\be F\|_{\vr,\nu(K)\times \nu(L)}^{m+n(|\al|+|\be|),\mu+|\al|+|\be|} \le C_1^{|\al|+|\be|} \|f\|_{\vr,K\times L}^{m,\mu}
\]
for all $\al,\be \in \N^n$. Hence for $\al,\be \in \N^n$ and $(u,v) \in \nu(K)\times \nu(L)$ we find
\begin{align*}
 |\p_u^\al &\p_v^\be F(u,v)| \\
&\le \|\p_u^\al \p_v^\be F\|_{\vr,\nu(K)\times \nu(L)}^{n(|\al|+|\be|),|\al|+|\be|} \vr^{n(|\al|+|\be|)} (n(|\al|+|\be|)+1)^{|\al|+|\be|} M_{n(|\al|+|\be|)}\\
&\le \|f\|_{\vr,K\times L}^{0,0} \, C_1^{|\al|+|\be|} \vr^{n(|\al|+|\be|)} (n(|\al|+|\be|)+1)^{|\al|+|\be|} M_{n(|\al|+|\be|)}\\
&\le C_2 \vr_1^{|\al|+|\be|} (|\al|+|\be|)! \, N_{|\al|+|\be|},
\end{align*}
for suitable constants $C_2$ and $\vr_1$. That implies $F \in C^N(\nu(\R^n)\times \nu(\R^n))$.
\qed\end{demo}

It was proved by Kostov \cite{Kostov89} that $\si(\R^n)$ is Whitney $1$-regular. 
Hence $\th(\R^{pn})=\si(\R^n) \times \cdots \times \si(\R^n)$ is Whitney $1$-regular as well.
It follows that, if $N$ is strongly regular, then $F$ can be extended to a function in 
$C^N(\R^{pn})$ (by Whitney's extension theorem; see \ref{CMjets}):

\begin{corollary*} 
Assume that $M$ is an increasing logarithmically convex sequences with $M_0=1$.
Let $N$ be a strongly regular DC-weight sequence.
For any  function $f \in C^M(\R^{pn})^{\prod_{j=1}^p \on{S}_n}$ there exists a function $F \in C^N(\R^{pn})$ 
such that $f = F \o \th$ if and only if 
\[
 \sup_{k \in \N_{>0}} \Big(\frac{M_{k n}}{N_k}\Big)^{\frac{1}{k}} < \infty.
\]
In particular: Any Gevrey function $f \in G^{1+\de}(\R^{pn})^{\prod_{j=1}^p \on{S}_n}$ (with $\de>0$) has the form $f = F \o \th$ with $F \in G^{1+\ga}(\R^{pn})$, 
where the exponent $\ga= \de n$ is minimal possible.
\end{corollary*}

\subsection{Invariant functions in Denjoy--Carleman classes} \label{ifDC}

\begin{theorem*}
Let $G$ be subgroup with finite order $m$ of $\on{GL}(V)$.  
Let $\si_1,\ldots,\si_p$ be a system of homogeneous generators of $\R[V]^G$ and put $\si=(\si_1,\ldots,\si_p) : V \to \R^p$.
Assume that $M$ and $N$ are DC-weight sequences.
Suppose that $N$ is strongly regular and that 
\begin{equation} \label{Gcond}
\sup_{k \in \N_{>0}} \Big(\frac{M_{k m}}{N_k}\Big)^{\frac{1}{k}} < \infty.
\end{equation}
Then for any $G$-invariant function $f \in C^M(V)^G$ there exists a function $F \in C^N(\R^p)$ 
such that $f = F \o \si$. 
In particular: Any $G$-invariant Gevrey function $f \in G^{1+\de}(V)^G$ (with $\de>0$) has the form $f = F \o \si$ 
with $F \in G^{1+\ga}(\R^p)$, where $\ga=\de m$.
\end{theorem*}

The proof of the theorem uses \ref{sym} and occupies the rest of the section. 
It is inspired by Barbançon and Raïs \cite{BR83} 
deploying Weyl's account \cite{Weyl39} of Noether's \cite{Noether16} proof of Hilbert's theorem.

\begin{remarks*}
(1) The condition \eqref{Gcond} implies that $C^M(U) \subseteq C^N(U)$ by \eqref{incl}. If additionally 
$\lim_{k\to \infty} (M_k/N_k)^{1/k} = 0$ then $C^M(U) \ne C^N(U)$, so there is a real loss of regularity. 

(2) The loss of regularity announced in the theorem is not minimal. For a particular group $G$, one may 
find much better Denjoy--Carleman regularity for $F$.

(3) The result is independent of the choice of generators $\si_i$, since any two choices differ by a 
polynomial diffeomorphism and the involved Denjoy--Carleman classes are stable under composition.
\end{remarks*}

\subsection{Reduction to the symmetric case} \label{reduct}

Let $V$ be a real vector space of finite dimension $n$ and 
let $G$ be a subgroup with finite order $m$ of $\on{GL}(V)$. The symmetric group $\on{S}_m$ acts in 
a natural way on $G$ by permuting the elements. This induces an action of $\on{S}_m$ on the 
space $F(G,\R)$ of functions defined in $G$ with values in $\R$ (for $\si \in \on{S}_m$ and $f \in F(G,\R)$
we have $\si.f =f \o \si^{-1}$). It can be identified with the standard representation $\rh$ of 
$\on{S}_m$ in $\R^m$. We obtain a natural action of $\on{S}_m$ on $E=F(G,\R) \otimes V$, the vector space 
of functions defined in $G$ with values in $V$. The corresponding representation $\pi$ is given by $\pi = n \rh$.

Let $L : V \to E$ be the linear injective mapping defined by $L(v) : g \mapsto g.v$ for $v \in V$.
We consider the pullback $L^* : F(E,\R) \to F(V,\R)$ (where $F(X,Y)$ denotes the space of functions defined in 
$X$ with values in $Y$). It is linear and maps $\on{S}_m$-invariant functions to $G$-invariant functions.
Hence it drops to a mapping $L^* : F(E,\R)^{\on{S}_m} \to F(V,\R)^G$. We define a linear mapping 
$J : F(V,\R) \to F(E,\R)$ by putting 
\[
J(f)(h) = \frac{1}{m} \sum_{g \in G} f(h(g))
\] 
for $f \in F(V,\R)$ and $h \in E = F(G,\R) \otimes V$. 
If we denote by $\on{ev}_g : E \to V$ the evaluation at $g \in G$, i.e., $\on{ev}_g(h)=h(g)$ for $h \in E$, 
then $J(f) = \frac{1}{m} \sum_{g \in G} \on{ev}_g^* f$. Thus, $J$ maps polynomials on $V$ to polynomials on $E$.
It is easy to check that $L^* \o J|_{F(V,\R)^G}= \on{id}$, so $J|_{F(V,\R)^G}$ is a section for 
$L^* : F(E,\R)^{\on{S}_m} \to F(V,\R)^G$.

Let $M$ be a DC-weight sequence.
It is easily seen that $L^*$ and $J$ are both continuous as mappings $L^* : C^M(E)^{\on{S}_m} \to C^M(V)^G$ 
and $J : C^M(V)^G \to C^M(E)^{\on{S}_m}$. 

Let $(\ta_1,\ldots,\ta_p)$ be a system of generators of the algebra $\R[E]^{\on{S}_m}$. 
Let $f \in C^M(V)^G$. If theorem \ref{ifDC} holds for $\pi$, there exists $F \in C^N(\R^p)$ 
(with suitable strongly regular DC-weight sequence $N$, see \ref{pi}) 
such that 
\[
J(f)(h)=F(\ta_1(h),\ldots,\ta_p(h))
\]
for all $h \in E$. Then
\[
f(v)=J(f)(L(v))=F(\si_1(v),\ldots,\si_p(v))
\]
for all $v \in V$, where $\si_i= L^*\ta_i$ for $1 \le i \le p$. It is clear from the above that 
the $\si_i= L^*\ta_i$ generate $\R[V]^G$.
This shows theorem \ref{ifDC} 
under the assumption that it holds for the representation $\pi$ (with suitable $N$).

\subsection{\label{refl}\!\!}

Let $W \subseteq \on{GL}(V)$ be a finite reflection group. Let $H$ be a $W$-invariant graded linear subspace 
of $\R[V]$ which is complementary to the ideal generated by the $W$-invariant polynomials with strictly 
positive degree. The bilinear mapping $(h,f) \mapsto h f$ induces an isomorphism of $W$-modules 
$H \otimes \R[V]^W \to \R[V]$ (see \cite[Ch. 5, 5.2, Thm. 2]{Bourbaki68}). So $\R[V]$ is a free $\R[V]^W$-module of rank $|W|$.  
 
Choose a basis $h_1,\ldots,h_{|W|}$ of $H$ consisting of homogeneous elements.
Let $w_1,\ldots,w_{|W|}$ denote the elements of $W$ (in some ordering). 
Since $\R[V] = H \R[V]^W$, we find that, for each $v \in V$, the cardinality of the orbit $W.v$ equals 
the rank of the matrix $(h_j(w_i.v))_{i,j}$. Since there are $v \in V$ with $|W.v|=|W|$, the polynomial
\[
\De(v):=\det (h_j(w_i.v))_{i,j}
\]
is not $0 \in \R[V]$.

\begin{lemma*}
Let $W=\on{S}_{m_1} \times \cdots \times \on{S}_{m_n}$ act in 
$V=\R^{m_1} \oplus \cdots \oplus \R^{m_n}$ by permuting the coordinates.
Then, for $v=(x_{1,1},\ldots,x_{1,m_1},\ldots,x_{n,1},\ldots,x_{n,m_n})$, we have 
\begin{equation} \label{Dpoly}
\De(v) = c \prod_{i=1}^n \prod_{1 \le j_i < k_i \le m_i} (x_{i,j_i}-x_{i,k_i})^{p_{i,j_i,k_i}}
\end{equation}
for some non-zero constant $c$ and positive integers $p_{i,j_i,k_i}$.
\end{lemma*}

\begin{demo}{Proof}
By definition, $\De(v)=0$ if and only if $v$ belongs to some reflecting hyperplane of $W$. 
It follows that each of the linear forms 
\begin{equation} \label{cL}
\cL:=\{x_{i,j_i}-x_{i,k_i} : 1 \le i \le n,1 \le j_i < k_i \le m_i\}
\end{equation}
divides $\De$. Since they are relatively prime, their product divides $\De$.
Suppose, for contradiction, there is a non-constant polynomial $P$ which is relatively prime with any of 
the linear forms in $\cL$ and divides $\De$. Without loss we switch to the complexification of the $W$-module $V$.
By Hilbert's Nullstellensatz, there is a positive integer $r$ such that $(\prod_{l \in \cL} l)^r$ belongs to 
the ideal generated by $\De$, a contradiction.
Hence the assertion.
\qed\end{demo}

\begin{remark*}
Actually, more is true: For any finite reflection group $W \subseteq \on{GL}(V)$ we have $\De=cJ^{|W|/2}$, 
where $c$ is a non-zero constant and $J=\prod_{l \in \cL_W} l$ with $\cL_W$ the set of linear forms with kernel a
reflection hyperplane of $W$. See \cite[4.2 + Appendix]{BR83}.
For us the above lemma will suffice.
\end{remark*}

\subsection{\label{modCM}\!\!}
Let $H$ and $h_1,\ldots,h_{|W|}$ be as in \ref{refl}.
The following proposition is a modification of \cite[3.3]{BR83}.

\begin{proposition*}
Let $M$ be a DC-weight sequence.
Let $W=\on{S}_{m_1} \times \cdots \times \on{S}_{m_n}$ act in 
$V=\R^{m_1} \oplus \cdots \oplus \R^{m_n}$ by permuting the coordinates.
Then $h_1,\ldots,h_{|W|}$ constitutes a basis of $C^M(V)$ considered as $C^M(V)^W$-module. 
\end{proposition*}

\begin{demo}{Proof}
Let $f \in C^M(V)$. There exists a sequence $(P_k)$ of polynomials which converges to $f$ in $C^M(V)$ 
(by \ref{dense}).
Since $h_1,\ldots,h_{|W|}$ is a basis of $\R[V]$ as $\R[V]^W$-module,
we can write $P_k = \sum_j h_j P_{k,j}$ with $P_{k,j} \in \R[V]^W$. For each $v \in V$, we obtain a system of 
$|W|$ equations
\[
P_k(w_i.v) = \sum_j h_j(w_i.v) P_{k,j}(v) \quad (1 \le i \le |W|).
\] 
Cramer's rule implies
\[
\De(v) P_{k,j}(v) = \sum_i \De_{ij}(v) P_k(w_i.v) \quad (1 \le j \le |W|),
\]
where the $\De_{ij}$ denote the cofactors of the matrix $(h_j(w_i.v))_{i,j}$.
The right-hand side of the single equations converges in $C^M(V)$ to the function 
\[
v \mapsto \sum_i \De_{ij}(v) f(w_i.v) \quad (1 \le j \le |W|),
\]
respectively (a straightforward computation shows that multiplication by a polynomial is continuous). 
Hence, each sequence $(\De P_{k,j})_k$ converges in $C^M(V)$. 
By proposition \ref{clideal} and lemma \ref{refl}, the ideal $\De C^M(V)$ generated by $\De$ is closed in $C^M(V)$.
Thus, there exist unique functions $f_j \in C^M(V)$ such that, for each $v$ and each $j$,
\begin{equation} \label{Dfj}
\De(v) f_j(v) = \sum_i \De_{ij}(v) f(w_i.v).
\end{equation} 
The $f_j$ are $W$-invariant: For each $w \in W$ there is $\ep_w \in \{0,1\}$ such that $\De(w.v)=(-1)^{\ep_w}\De(v)$ for all $v \in V$.
Since the polynomials $P_{k,j}$ are $W$-invariant and evaluation at points is continuous, we find 
\[
 (-1)^{\ep_w}\De(v) f_j(w.v) = (-1)^{\ep_w}\De(v) f_j(v)
\]
and thus
\[
f_j(w.v)=f_j(v) 
\]
on the open dense subset $\{v : \De(v) \ne 0\}$, and hence everywhere.
From \eqref{Dfj} we obtain
\[
f(v) = \sum_j h_j(v) f_j(v)
\]
on the open dense subset $\{v : \De(v) \ne 0\}$, and hence everywhere.
\qed\end{demo}

\begin{remark*}
Using remark \ref{refl}, we find that this proposition is true for any finite reflection group $W \subseteq \on{GL}(V)$.
\end{remark*}

\subsection{Theorem \ref{ifDC} for the representation $\pi : \on{S}_m \to \on{GL}(\R^{nm})$}  \label{pi}
Let $G$ be a subgroup of $W=\on{S}_{m_1} \times \cdots \times \on{S}_{m_n}$ acting in 
$V=\R^{m_1} \oplus \cdots \oplus \R^{m_n}$ by permuting the coordinates. 
Let $H$ be (as in \ref{refl}) a $W$-invariant graded linear subspace 
of $\R[V]$ which is complementary to the ideal generated by the $W$-invariant polynomials with strictly 
positive degree. 
Consider a basis $(h_1,\ldots,h_r)$ of $H^G$. 
By proposition \ref{modCM}, we find that $(h_1,\ldots,h_r)$ constitutes a basis of 
$C^M(V)^G$ considered as $C^M(V)^W$-module. 

By the reduction in \ref{reduct}, in order to prove theorem \ref{ifDC} it suffices to consider the representation $\pi : \on{S}_m \to \on{GL}(\R^{nm})$.
Let $\ta_1,\ldots.\ta_p$ and $\th_1,\ldots,\th_{nm}$ be systems of homogeneous generators of $\R[\R^{nm}]^{\on{S}_m}$ and
$\R[\bigoplus_{j=1}^n \R^m]^{\prod_{j=1}^n \on{S}_m}$, respectively, and consider $\ta=(\ta_1,\ldots,\ta_p) : \R^{nm} \to \R^p$ and 
$\th=(\th_1,\ldots,\th_{nm}) : \R^{nm} \to \R^{nm}$. By the previous paragraph and corollary \ref{sym}, 
each $f \in C^M(\R^{nm})^{\on{S}_m}$ has the form
\[
 f=\sum_{j=1}^r h_j f_j = \sum_{j=1}^r (H_j \o \ta)(F_j \o \Th \o \ta),
\]
where $h_j \in \R[\R^{nm}]^{\on{S}_m}$, $f_j \in C^M(\R^{nm})^{\prod_{j=1}^n \on{S}_m}$,
$H_j \in \R[\R^p]$, $F_j \in C^N(\R^{nm})$, and $\Th$ is the polynomial mapping given by $\th=\Th \o \ta$.
Note that $N$ is a strongly regular DC-weight sequence satisfying
\[
 \sup_{k \in \N_{>0}} \Big(\frac{M_{k m}}{N_k}\Big)^{\frac{1}{k}} < \infty.
\]
This completes the proof of theorem \ref{ifDC}.

\section{Equivariant mappings in Denjoy--Carleman classes} \label{secequi}

We give an application of theorem \ref{ifDC} to the representation of 
equivariant mappings in Denjoy--Carleman classes. We follow standard techniques.

\subsection{\label{emDC}\!\!}
Let $V_1$ and $V_2$ be real finite dimensional representations of a compact Lie group $G$. 
It is well-known that the set $\on{Pol}(V_1,V_2)^G$ of $G$-equivariant polynomial mappings from $V_1$ to 
$V_2$ is finitely generated as module over $\R[V_1]^G$. 

Let $M$ be a DC-weight sequence.
We denote by $C^M(V_1,V_2)^G$ the set of $G$-equivariant $C^M$-mappings $f : V_1 \to V_2$.

\begin{theorem*}
Let $V_1$ and $V_2$ be representations of a finite group $G$ with order $m$.  
Let $\si_1,\ldots,\si_p$ be a system of homogeneous generators of $\R[V_1]^G$ and put $\si=(\si_1,\ldots,\si_p)$.
Let $P_1,\ldots,P_l$ be a system of generators of the $\R[V_1]^G$-module $\on{Pol}(V_1,V_2)^G$.
Assume that $M$ and $N$ are DC-weight sequences. 
Suppose that $N$ is strongly regular and that 
\[
 \sup_{k \in \N_{>0}} \Big(\frac{M_{k m}}{N_k}\Big)^{\frac{1}{k}} < \infty.
\]
Then for each $f \in C^M(V_1,V_2)^G$ there exists an $L(f) \in (C^N(\R^p))^l$
such that $f = \sum_{j=1}^l (L(f)_j \o \si) P_j$.
\end{theorem*}

\begin{demo}{Proof}
The dual $V_2^*$ of $V_2$ carries the dual $G$-action given by $g.l=l \o g^{-1}$. 
Let $f \in C^M(V_1,V_2)^G$ and consider the $G$-invariant function $H_f : V_1 \times V_2^* \to \R$
given by $H_f(v,l)=l(f(v))$. 
So $H_f \in C^M(V_1 \times V_2^*)^G$ and, by theorem \ref{ifDC}, 
there exists $L_f \in C^N(\R^q)$ such that $H_f=L_f \o \ta$, where $\ta=(\ta_1,\dots,\ta_q)$ and 
$\ta_1,\dots,\ta_q$ generate $\R[V_1 \times V_2^*]^G$. Taking the derivative with respect to the second 
component gives 
\[
f(v) = \sum_{i=1}^q \p_i L_f(\ta(v,0)) d_2 \ta_i(v,0).
\]
Since $v \mapsto d_2 \ta_i(v,0)$ is a $G$-equivariant polynomial mapping, there exist $h_{ij} \in \R[V_1]^G$ 
such that $d_2 \ta_i(v,0) = \sum_{j=1}^l h_{ij}(v) P_j(v)$. Since $v \mapsto \ta(v,0)$ is $G$-invariant, 
there is a polynomial mapping $\th : \R^p \to \R^q$ with $\ta(v,0)= \th(\si(v))$.
Then
\[
L(f) := \Big(\sum_{i=1}^q (\p_i L_f \o \th) h_{ij}\Big)_{1 \le j \le l}
\] 
has the required properties.
\qed\end{demo}

\section{Polar representations} \label{secpolar}

\subsection{Polar representations} \label{polar}
\cite{DK85}, \cite{PT88}, \cite{Terng85}
A real finite dimensional orthogonal representation $\rh : G \rightarrow \on{O}(V)$ of a Lie group $G$ 
is called \emph{polar}, 
if there exists a linear subspace $\Si \subseteq V$, 
called a \emph{section}, 
which meets each orbit orthogonally. 
The trace of the $G$-action in $\Si$ is the action of the \emph{generalized Weyl group} 
$W(\Si) = N_G(\Si)/Z_G(\Si)$, where 
$N_G(\Si) := \{g \in G : \rh(g)(\Si) = \Si\}$ and 
$Z_G(\Si) := \{g \in G : \rh(g)(s) = s ~\text{for all}~ s \in \Si\}$. 
The generalized Weyl group is a finite group. 
If $\Si'$ is a different section, 
then there is an isomorphism $W(\Si) \to W(\Si')$ 
induced by an inner automorphism of $G$.

The following generalization of Chevalley's restriction theorem 
is due to Dadok and Kac \cite{DK85} 
and independently to Terng \cite{Terng85}.

\begin{theorem*} 
Assume that $G$ is a compact Lie group.
Then restriction induces an isomorphism of algebras between $\R[V]^G$  
and $\R[\Si]^{W(\Si)}$. 
\end{theorem*}

\subsection{Invariant functions in Denjoy--Carleman classes} \label{polDC}
We generalize theorem \ref{ifDC} to polar representations..

\begin{theorem*}
Let $G \to \on{O}(V)$ be a polar representation of a compact Lie group $G$, with section $\Si$ and 
generalized Weyl group $W=W(\Si)$. Write $m=|W|$.  
Let $\si_1,\ldots,\si_p$ be a system of homogeneous generators of $\R[V]^G$ and put $\si=(\si_1,\ldots,\si_p)$.
Assume that $M$ and $N$ are DC-weight sequences.
Suppose that $N$ is strongly regular and that 
\[
 \sup_{k \in \N_{>0}} \Big(\frac{M_{k m}}{N_k}\Big)^{\frac{1}{k}} < \infty.
\]
Then for any $G$-invariant function $f \in C^M(V)^{G}$ there exists a function $F \in C^N(\R^p)$ 
such that $f = F \o \si$. 
In particular: Any $G$-invariant Gevrey function $f \in G^{1+\de}(V)^G$ (with $\de>0$) has the form $f = F \o \si$ 
with $F \in G^{1+\ga}(\R^p)$, where $\ga=\de m$.
\end{theorem*}

\begin{demo}{Proof}
Let $f \in C^M(V)^{G}$.  
By theorem \ref{polar}, the restrictions $\si_1|_\Si,\ldots,\si_p|_\Si$ generate $\R[\Si]^W$ and 
$\si(V)=\si|_\Si(\Si)$. Since $f|_\Si \in C^M(\Si)^W$, theorem \ref{ifDC} implies that there is a 
$F \in C^N(\R^p)$ such that $f|_\Si = F \o \si|_\Si$, and, hence, $f = F \o \si$. 
\qed\end{demo}

\subsection{\label{invext}\!\!} 
In the situation of \ref{polDC} we have:

\begin{theorem*}
Each $f \in C^M(\Si)^W$ (resp.\ $G^{1+\de}(\Si)^W$) has an extension in $C^N(V)^G$ (resp.\ $G^{1+\ga}(V)^G$). 
\end{theorem*}

\begin{demo}{Proof}
Let $f \in C^M(\Si)^W$. Choose a system of homogeneous generators $\ta_1,\ldots,\ta_p$ of $\R[\Si]^W$. 
By theorem \ref{ifDC}, there is an $F \in C^N(\R^p)$ such that $f = F \o (\ta_1,\ldots,\ta_p)$. 
Each $\ta_i$ extends to a polynomial $\tilde \ta_i \in \R[V]^G$, by theorem \ref{polar}. 
So $\tilde f := F \o (\tilde \ta_1,\ldots,\tilde \ta_p)$ is a $G$-invariant extension of $f$ 
belonging to $C^N(V)$.
\qed\end{demo}

\subsection{Basic differential forms in Denjoy--Carleman classes} \label{bdfpolar}
Let $G \to \on{O}(V)$ be a polar representation of a compact Lie group $G$, with section $\Si$ and 
generalized Weyl group $W=W(\Si)$.
A differential form $\om \in \Om^p(V)$ is called \emph{$G$-invariant} if $(l_g)^* \om = \om$ for all $g \in G$, where $l_g(x)=g.x$,
and \emph{horizontal} if it kills each vector tangent to a $G$-orbit, i.e., $i_{\ze_X} \om =0$ for all 
$X \in \mathfrak g := \on{Lie}(G)$, where $\ze$ is the fundamental vector field mapping ($\ze_X(x) = T_e (l^x).X$ with $l^x(g)=g.x$).
Denote by $\Om_{\on{hor}}^p(V)^G$ the space of all horizontal $G$-invariant $p$-forms on $V$. Its elements 
are also called \emph{basic $p$-forms}.

It is proved in \cite{Michor96B, Michor97B} that the restriction of differential forms induces an 
isomorphism between $\Om_{\on{hor}}^p(V)^G$ and $\Om^p(\Si)^W$.

Let $M$ be a DC-weight sequence. We may consider $p$-forms $\om$ on $V$ of 
Denjoy--Carleman class $C^M$. Let us denote the space of such forms $\om$ by $(\Om^M)^p(V)$.
A careful inspection of the proofs in \cite{Michor96B, Michor97B} shows that we can deduce the 
following theorem in an analog manner: 
\begin{enumerate}
\item[(i)] The statement in \cite[3.2]{Michor96B} is true in Denjoy--Carleman classes $C^M$ as well:
{\it Let $l \in V^*$ and let $f \in C^M(V)$ with $f|_{l^{-1}(0)}=0$. 
Then there exists a unique $h \in C^M(V)$ such that $f = l \cdot h$.}
See the proof of proposition \ref{clideal}.
\item[(ii)] In \cite[3.7]{Michor96B} instead of Schwarz's theorem we use theorem \ref{ifDC}.
\end{enumerate}
The rest works without change and yields:

\begin{theorem*}
Let $G \to \on{O}(V)$ be a polar representation of a compact Lie group $G$, with section $\Si$ and 
generalized Weyl group $W=W(\Si)$. Put $m=|W|$.
Assume that $M$ and $N$ are DC-weight sequences.
Suppose that $N$ is strongly regular and that 
\[
 \sup_{k \in \N_{>0}} \Big(\frac{M_{k m}}{N_k}\Big)^{\frac{1}{k}} < \infty.
\]
Then each $\om \in (\Om^M)^p(\Si)^W$ has an extension in $(\Om^N)_{\on{hor}}^p(V)^G$. \qed
\end{theorem*}

\begin{remark*}
Obviously, restriction of differential forms does in general not map forms in $(\Om^N)_{\on{hor}}^p(V)^G$ 
to forms in $(\Om^M)^p(\Si)^W$. So we cannot expect to obtain an isomorphism as in the smooth case.
\end{remark*}

\section{Proper $G$-manifolds with sections} \label{secGmf}

In this section $X$ always denotes a connected complete Riemannian $G$-manifold, with 
effective and isometric $G$-action.

\subsection{Sections}
\cite{PT88}
Let $X$ be a proper Riemannian $G$-manifold. A connected closed submanifold $\Si$ of $X$ is 
called a \emph{section} for the $G$-action, if it meets all $G$-orbits orthogonally. Each section is a 
totally geodesic submanifold. Analogously with \ref{polar} we define the \emph{generalized Weyl group} 
$W(\Si):= N_G(\Si)/Z_G(\Si)$ which turns out to be a discrete group acting properly on $\Si$. 
If $\Si'$ is a different section, then there is an isomorphism $W(\Si) \to W(\Si')$ induced by an inner 
automorphism of $G$.

\subsection{Invariant functions in Denjoy--Carleman classes}
In the smooth case, restriction induces an isomorphism $C^\infty(X)^G \cong C^\infty(\Si)^{W(\Si)}$, 
by \cite{PT87}. We show an analog result in Denjoy--Carleman classes. From now on all manifolds are real analytic.

\begin{theorem*}
Let $X$ be a real analytic proper Riemannian $G$-manifold with section $\Si$ 
and Weyl group $W=W(\Si)$.
Suppose that 
\[
 m:=\sup_{x \in \Si} |W_x| < \infty. 
\]
Assume that $M$ and $N$ are DC-weight sequences.
Suppose that $N$ is strongly regular and that 
\[
 \sup_{k \in \N_{>0}} \Big(\frac{M_{k m}}{N_k}\Big)^{\frac{1}{k}} < \infty.
\]
Then each $f \in C^M(\Si)^W$ (resp.\ $G^{1+\de}(\Si)^W$) has an extension in $C^N(X)^G$ (resp.\ $G^{1+\de m}(X)^G$).
\end{theorem*}

\begin{demo}{Proof}
Let $f \in C^M(\Si)^W$. It is well-known (e.g.\ \cite{PT88}) that each $W$-invariant 
continuous (smooth) function in $\Si$ has a unique continuous (smooth) $G$-invariant extension. 
Let $\tilde f$ be the extension of $f$. 
We show that $\tilde f$ represents an element in $C^N(X)^G$.  
Let $x \in X$. Without loss we may assume that $x \in \Si$ (since the action is real analytic). 
Let $S_x$ be a normal slice at $x$. 
Then, by the slice theorem, $G.S_x$ and $G \times_{G_x} S_x$ are real analytically $G$-isomorphic and 
$G \times S_x \to G \times_{G_x} S_x$ is a real analytic surjective submersion. 
Thus, it suffices to show that $\tilde f|_{S_x}$ 
belongs to $C^N(S_x)$. We can choose a ball $B \subseteq T_x S_x$ around $0_x$ such that 
$B \cong S_x$ and $T_x \Si \cap B \cong \Si \cap S_x$. Then the $G_x$-action on $S_x$ is up to a real 
analytic isomorphism a polar representation with section $T_x \Si$ and Weyl group $W_x$ (e.g.\ \cite{PT88}). 
So the assertion follows from theorem \ref{invext}.
\qed\end{demo}

\subsection{Basic differential forms in Denjoy--Carleman classes}
In the smooth case, the restriction of differential forms induces an 
isomorphism between $\Om_{\on{hor}}^p(X)^G$ and $\Om^p(\Si)^W$, by \cite{Michor96B, Michor97B}. This is derived from the analog result 
for polar representations with the help of the slice theorem.

Let $X$ be a real analytic proper Riemannian $G$-manifold with section $\Si$ and Weyl group $W=W(\Si)$.
Suppose that 
\[
 m:=\sup_{x \in \Si} |W_x| < \infty. 
\]
Let $M$ and $N$ be a DC-weight sequences.
Suppose that $N$ is strongly regular and that
\[
 \sup_{k \in \N_{>0}} \Big(\frac{M_{k m}}{N_k}\Big)^{\frac{1}{k}} < \infty.
\]

It turns out that we are able to apply the same arguments as in \cite{Michor96B, Michor97B} in order to deduce a similar generalized statement 
from theorem \ref{bdfpolar}.  
All mappings occurring while 
applying the slice theorem in \cite[part 4]{Michor96B} are real analytic and may, therefore, be taken over without change.
Hence we may reduce to the slice representations $G_x \to \on{O}(T_x S_x)$ which are polar with Weyl group $W_x$ and we may apply theorem \ref{bdfpolar}.

Following the final step of the proof \cite[4.2]{Michor96B} we glue local differential forms $\om^{x_n} \in (\Om^N)_{\on{hor}}^p(G.S_{x_n})^G$ to 
a form $\tilde \om \in (\Om^N)_{\on{hor}}^p(X)^G$.  
This is done, using a method of Palais \cite[4.3.1]{Palais61}, by constructing a suitable partition of 
unity consisting of $G$-invariant functions.
More precisely:
There exists a sequence $(x_n)_{n\in \N}$ of points in $\Si$ and open neighborhoods of $x_n$ in $\Si$ whose projections form a locally finite open covering of the 
orbit space $X/G \cong \Si/W$,
and there exists a partition of unity $f_n$ consisting of $G$-invariant functions with $\on{supp}(f_n) \subseteq G.S_{x_n}$.
The construction of the $f_n$ is as follows: There exist neighborhoods $x_n \in K_n$ with compact closure in $S_{x_n}$ such that their projection forms a covering of $X/G$. Let $f_n$ be a non-negative function on $S_{x_n}$ positive on $K_n$ and with compact support in $S_{x_n}$. By averaging we may assume that $f_n$ is $G_{x_n}$-invariant. Define $f_n(g.s)=f_n(s)$ for $g \in G$ and $s \in S_{x_n}$ and $f_n(x) =0$ for $x \not\in G.S_{x_n}$.
Since there are $C^N$ partitions of unity (by \ref{2.3}) and since averaging over the slice representation $G_{x_n} \to \on{O}(T_{x_n} S_{x_n})$ 
(which is $G_{x_n}$-equivariantly real analytically isomorphic to the $G_{x_n}$-manifold $S_{x_n}$)
preserves the Denjoy--Carleman class (by lemma \ref{average} below), the functions $f_n$ can be chosen in $C^N(X)^G$. 
Thus $\tilde \om = \sum_n f_n \om^{x_n}  \in (\Om^N)_{\on{hor}}^p(X)^G$, and we obtain:

\begin{theorem*}
Let $X$ be a real analytic proper Riemannian $G$-manifold with section $\Si$ and Weyl group $W=W(\Si)$.
Suppose that 
\[
 m:=\sup_{x \in \Si} |W_x| < \infty. 
\]
Let $M$ and $N$ be a DC-weight sequences.
Suppose that $N$ is strongly regular and that
\[
 \sup_{k \in \N_{>0}} \Big(\frac{M_{k m}}{N_k}\Big)^{\frac{1}{k}} < \infty.
\]
Then each $\om \in (\Om^M)^p(\Si)^W$ has an extension in $(\Om^N)_{\on{hor}}^p(X)^G$. \qed
\end{theorem*}

\begin{lemma}\label{average}
Let $G \to \on{O}(V)$ be a real finite dimensional representation of a compact Lie group $G$. Let $M$ be a DC-weight sequence. If $f \in C^M(V)$ then 
\[
\tilde f(x)=\int_G f(g.x) dg
\] 
(where $dg$ denotes Haar measure) belongs to $C^M(V)^G$.
\end{lemma}

\begin{demo}{Proof}
We write $l_g : V \to V, x \mapsto g.x$ for the linear action of $g \in G$. By choosing a basis we identify $V=\R^n$.
Let $K \subseteq V$ be compact. It suffices to show that for each positive $\vr = \vr(f,G.K)$ there exists a positive $\bar \vr$ such that 
\begin{equation}\label{eqaverage}
 \|f \o l_g\|_{\bar \vr,K} \le \|f\|_{\vr,G.K}
\end{equation}
for all $g \in G$. By Fa\`a di Bruno (\cite{FaadiBruno1855} for the 1-dimensional version)
\[
\frac{\p^\ga (f \o l_g)(x)}{\ga!} 
= \sum_{\substack{\be_i \in \N^n \backslash \{0\}\\ \al=\be_1+\cdots+\be_n\\ \ga=(|\be_1|,\ldots,|\be_n|)}} 
\frac{1}{\be_1! \cdots \be_n!} \, \p^\al f(g.x) \, (\p_1 l_g(x))^{\be_1} \cdots (\p_n l_g(x))^{\be_n},
\]
where $\p_i l_g(x) = (\p_i (l_g)_1(x),\dots,\p_i (l_g)_n(x))$. So we find
\[
\frac{|\p^\ga (f \o l_g)(x)|}{|\ga|!\, M_{|\ga|}} 
\le \sum
\frac{|\al|!}{\be_1! \cdots \be_n!} \,  \frac{|\p^\al f(g.x)|}{|\al|!\, M_{|\al|}} \, \|l_g\|^{|\al|}, 
\]
where $\|l_g\|$ denotes the operator norm of $l_g$. Put
\[
\mu := \max_{g \in G} \|l_g\|.
\]
Then we obtain \eqref{eqaverage} by defining
\[
\bar \vr := n^2 \mu \vr. 
\] 
This completes the proof.
\qed\end{demo}

\end{document}